\documentclass[english]{amsart}
\usepackage[T1]{fontenc}
\usepackage{textcomp}
\usepackage{amsthm}
\usepackage{amstext}
\usepackage{amssymb}
\usepackage{graphicx}

\makeatletter

\newcommand{\lyxmathsym}[1]{\ifmmode\begingroup\def\b@ld{bold}
  \text{\ifx\math@version\b@ld\bfseries\fi#1}\endgroup\else#1\fi}

\providecommand{\tabularnewline}{\\}

\numberwithin{equation}{section}
\numberwithin{figure}{section}

\makeatother

\usepackage{babel}
\begin{document}

\title{Conjecture on the value of $\pi\left(10^{26}\right)$, the number
of primes $<10^{26}$}

\author{Vladimir Pletser}

\address{European Space Research and Technology Centre \\
ESA-ESTEC P.O. Box 299 \\
NL-2200 AG Noordwijk \\
The Netherlands}

\address{European Space Research and Technology Centre \\
ESA-ESTEC P.O. Box 299 \\
NL-2200 AG Noordwijk \\
The Netherlands}

\email{Vladimir.Pletser@esa.int}


\begin{abstract}
Based on the first 25 known values of $\pi\left(10^{n}\right)$, the
number of primes less than $10^{n}$ with $1\leq n\leq25$, we propose
a conjectured value range of $\pi\left(10^{26}\right)$ calculated
by using polynomial interpolations with two corrective functions obtained
by Thiele interpolations on relative differences of exact and interpolated
values of $\pi\left(10^{n}\right)$. The conjectured range value is
in agreement with values obtained by the Eulerian logarithmic integral
and with the Riemann functions.
\end{abstract}
\maketitle

\section{Introduction}

For $x\in\mathbb{Z}^{+}$, it is known (see e.g. \cite{1}) that if
$\pi(x)$ is the number of primes$\leq x$, then 
\begin{equation}
lim_{x\rightarrow\infty}\left(\frac{\pi(x)}{\frac{x}{log\left(x\right)}}\right)=1\label{eq:1.1}
\end{equation}
meaning in the general sense that the relative error of approximating
$\pi(x)$ by $\left(\frac{x}{log\left(x\right)}\right)$ approaches
$0$ as $x$ approaches infinity. Several better approximations to
$\pi(x)$ for $x\lesssim10^{25}$ are given for example by the offset
logarithmic integral or Eulerian logarithmic integral (see e.g. \cite{2})
\begin{equation}
Li(x)=li\left(x\right)-li\left(2\right)=\intop_{2}^{x}\frac{dt}{log\left(t\right)}\label{eq:1.2}
\end{equation}
in function of the logarithmic integral $li(x)=\intop_{0}^{x}\frac{dt}{log\left(t\right)}$,
or by the Riemann function (see e.g. \cite{3})
\begin{equation}
R\left(x\right)=\sum_{j=1}^{\infty}\frac{\mu\left(j\right)}{j}li\left(x^{\frac{1}{j}}\right)\label{eq:1.3}
\end{equation}
where $\mu\left(j\right)$ is the Moebius function and $j\in\mathbb{Z}^{+}$,
or by an even better function
\begin{equation}
R\left(x\right)-\frac{1}{log\left(x\right)}+\frac{1}{\pi}arctan\left(\frac{\pi}{log\left(x\right)}\right)\label{eq:1.4}
\end{equation}

Other approximation formulas to $\pi(x)$ include Legendre's \cite{4},
Lehmer's \cite{5}, and Meissel's \cite{6}.

Values of $\pi(10^{n})$ are known for $1\leq n\leq25$ (see \cite{7}).
The values of $\pi(10^{24})$ and $\pi(10^{25})$ were calculated
relatively recently, respectively in 2010 \cite{8} and in June 2013
\cite{9}. 

In this paper, we propose a simple method based on polynomial interpolations
of the known 25 first values of $\pi\left(10^{n}\right)$ introducing
two corrective functions to calculate a conjectured value of $\pi\left(10^{26}\right)$.

\section{Method}

The method we use includes six steps, some involving interpolations%
\footnote{Curve Fitting package of Maple 16.00 with 50 digits precision%
}.

Step 1: First, polynomial functions 
\begin{equation}
P_{n}\left(x\right)=\sum_{i=0}^{n-1}\left(a_{i,n}x^{i}\right)\label{eq:2.1}
\end{equation}

with coefficients $a_{i,n}\in\mathbb{R}$ and $i,n,x\in\mathbb{Z}^{+}$,
are calculated by polynomial interpolations over the values of $\pi(10^{n})$
up to $n$ for $2\leq n\leq25$, yielding 23 polynomials, given in
Appendix 1, such that the equality
\begin{equation}
\pi(10^{x})=P_{n}\left(x\right)\label{eq:2.2}
\end{equation}

holds exactly or $1\leq x\leq n$.

Step 2: For each value of $n$ until $n=24$, the next value for $x=n+1$
is calculated by (\ref{eq:2.1}) yielding $P_{n}\left(n+1\right)$
and compared to the value of $\pi(10^{n+1})$. The relative differences

\begin{equation}
\delta_{n+1}=\frac{\pi(10^{n+1})-P_{n}\left(n+1\right)}{\pi(10^{n+1})}\label{eq:2.3}
\end{equation}

are then calculated. As $P_{n}\left(n+1\right)<\pi(10^{n+1})$ for
$2\leq n\leq24$, $\delta_{n+1}$ is such that $0<\delta_{n+1}<1$
and decreases with increasing values of $n$ (see Fig. 2.1).

\begin{figure}
\begin{raggedright}
\caption{$\delta_{n+1}$in function of $\left(n+1\right)$}

\par\end{raggedright}

\raggedright{}\includegraphics[width=12cm]{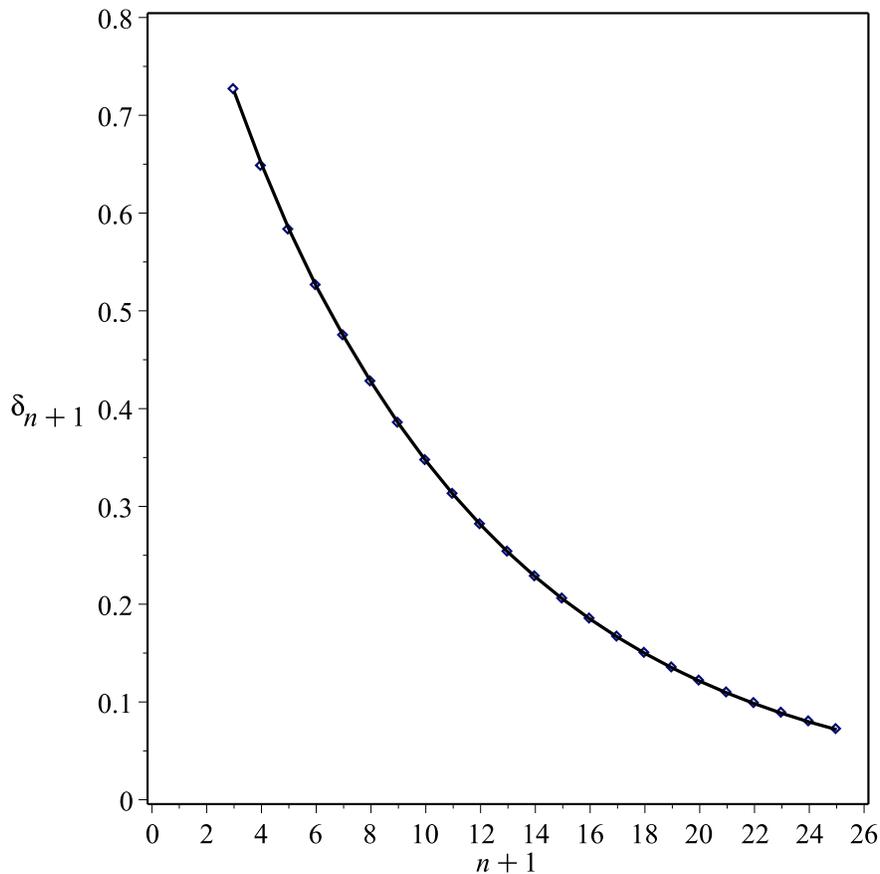}
\end{figure}

\begin{flushleft}
Step 3: For each set value of $n$, corrective functions $\Phi_{n+1}$
in the form of continued fractions 
\begin{equation}
\Phi_{n+1}\left(x\right)=c_{1}+\cfrac{x-3}{c_{2}+\cfrac{x-4}{c_{3}+...\underset{\ddots\,...+\cfrac{x-\left(n-3\right)}{c_{n-4}+\cfrac{x-\left(n-2\right)}{c_{n-3}+\cfrac{x-\left(n-1\right)}{c_{n-2}+Kx}}}}{}}}\label{eq:2.4}
\end{equation}

\par\end{flushleft}

with coefficients $c_{i},K\in\mathbb{R}$ and $i,n,x\in\mathbb{Z}^{+}$,
$1\leq i\leq n-2$, are calculated by Thiele interpolations over the
$\left(n-1\right)$ values of $\delta_{n+1}$ for $3\leq n\leq24$,
yielding 22 corrective functions $\Phi_{n+1}$. The coefficients $c_{i}$
and $K$ of $\Phi_{n+1}$ for $n=24$ are given in Appendix 2.

Step 4: The predictive power of using the 22 corrective functions
$\Phi_{n+1}$ to calculate a next interpolated value of $\pi(10^{x+1})$
knowing the first $x$ values of $\pi(10^{x})$ with $1\leq x\leq n$
for all values of $n$, $1\leq n\leq24$, is verified by taking the
nearest integer to $\left(P_{x}\left(x+1\right)/\left(1-\Phi_{n+1}\left(x+1\right)\right)\right)$
and compared to $\pi(10^{x+1})$. It yields that the equality
\begin{equation}
\pi(10^{x+1})=Round\left(\frac{P_{x}\left(x+1\right)}{1-\Phi_{n+1}\left(x+1\right)}\right)\label{eq:2.5}
\end{equation}

holds exactly for all values of $1\leq x\leq n$ and of $n$, $1\leq n\leq24$,
where $\Phi_{4}$ should be taken for $1\leq n\leq3$ and $Round(X)$
is the nearest integer to $X$ ($X\in\mathbb{R}$).

Step 5: To interpolate the next unknown value for $n=26$ of $\pi(10^{26})$,
as $\Phi_{26}$ is obviously unknown, we use a relation similar to
(\ref{eq:2.5}) with $x=25$ and $\Phi_{25}$ instead of $\Phi_{26}$,
to obtain a very approximate interpolated value of $\pi(10^{26})$
\begin{equation}
Round\left(\frac{P_{25}\left(26\right)}{1-\Phi_{25}\left(26\right)}\right)=1699246738822618041025224\label{eq:2.6}
\end{equation}

Step 6: To improve this interpolation, we have to find an additional
correction to the corrective function $\Phi_{25}$. We performed a
similar interpolation in May 2013 (unpublished) on the first 24 values
of $\pi(10^{n})$ before the announcement of the calculated value
of $\pi(10^{25})=176846309399143769411680$ by J. Buethe, J. Franke,
A. Jost and T. Kleinjung on 1st June 2013 \cite{9}. This previous
interpolation yielded an interpolated value of 
\begin{equation}
Round\left(\frac{P_{24}\left(25\right)}{1-\Phi_{24}\left(25\right)}\right)=176846307027334692763889<\pi(10^{25})\label{eq:2.7}
\end{equation}

giving a relative difference 
\begin{equation}
\delta_{n+1}^{\lyxmathsym{\textasciiacute}}=\frac{\pi(10^{n+1})-Round\left(\frac{P_{n}\left(n+1\right)}{1-\Phi_{n}\left(n+1\right)}\right)}{\pi(10^{n+1})}\label{eq:2.8}
\end{equation}

of $\delta_{25}^{\lyxmathsym{\textasciiacute}}\approx1.34117\times10^{-8}$.

To estimate the error made by using $\Phi_{n}$ instead of $\Phi_{n+1}$
in (\ref{eq:2.5}) when interpolating to find a value of $\pi(10^{n+1})$,
the relative differences $\delta_{n+1}^{\lyxmathsym{\textasciiacute}}$
are calculated for $4\leq n\leq24$, yielding 21 relative differences
$\delta_{n+1}^{\lyxmathsym{\textasciiacute}}$ (\ref{eq:2.8}), whose
approximate values are reported in Table 1. Their absolute value decreases
with increasing $n$ as shown in Fig. 2.2.

\begin{table}
\caption{Relative differences $\delta_{n+1}^{\lyxmathsym{\textasciiacute}}$}

\begin{tabular}{|c|c|c|c|c|c|}
\hline 
$n+1$ & $\delta_{n+1}^{\lyxmathsym{\textasciiacute}}$ & $n+1$ & $\delta_{n+1}^{\lyxmathsym{\textasciiacute}}$ & $n+1$ & $\delta_{n+1}^{\lyxmathsym{\textasciiacute}}$\tabularnewline
\hline 
\hline 
5 & $3.10676\times10^{-2}$ & 12 & $1.23708\times10^{-4}$ & 19 & $5.44808\times10^{-7}$\tabularnewline
\hline 
6 & $-4.15297\times10^{-3}$ & 13 & $-3.31760\times10^{-5}$ & 20 & $1.72692\times10^{-7}$\tabularnewline
\hline 
7 & $2.85895\times10^{-5}$ & 14 & $4.84531\times10^{-5}$ & 21 & $-3.64154\times10^{-7}$\tabularnewline
\hline 
8 & $7.32975\times10^{-4}$ & 15 & $-1.62333\times10^{-5}$ & 22 & $-7.56044\times10^{-8}$\tabularnewline
\hline 
9 & $5.14623\times10^{-3}$ & 16 & $1.11746\times10^{-6}$ & 23 & $-9.45864\times10^{-9}$\tabularnewline
\hline 
10 & $-8.04909\times10^{-3}$ & 17 & $-4.32482\times10^{-6}$ & 24 & $2.62139\times10^{-8}$\tabularnewline
\hline 
11 & $5.18791\times10^{-3}$ & 18 & $-4.83262\times10^{-6}$ & 25 & $1.34117\times10^{-8}$\tabularnewline
\hline 
\end{tabular}
\end{table}
\begin{figure}
\caption{$\left|\delta_{n+1}^{\lyxmathsym{\textasciiacute}}\right|$ in function
of $n+1$}

\includegraphics[width=12cm]{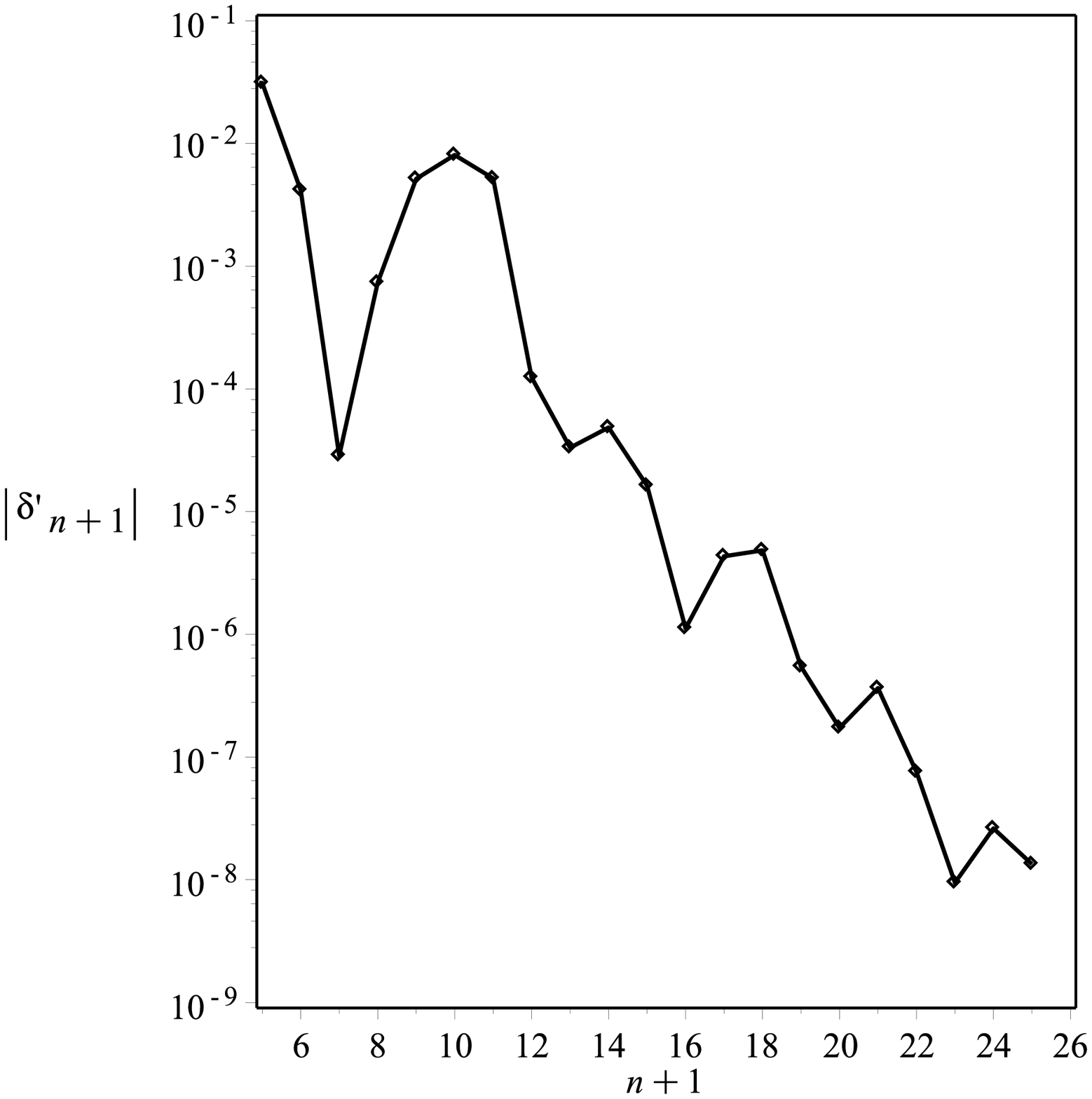}
\end{figure}

Assuming that$\left|\delta_{n+1}^{\lyxmathsym{\textasciiacute}}\right|$
continues to decrease with $n$, an interpolated value of $\left|\delta_{26}^{\lyxmathsym{\textasciiacute}}\right|\lesssim10^{-8}$
can be expected. A new Thiele interpolation over the last five values
of $\left|\delta_{n+1}^{\lyxmathsym{\textasciiacute}}\right|$ for
$20\leq n\leq24$ yields a second corrective function $\left|\psi_{n+1}\left(n\right)\right|$

\begin{flushleft}
\begin{equation}
\left|\psi_{n+1}\left(n\right)\right|=d_{1}+\cfrac{n-20}{d_{2}+\cfrac{n-21}{d_{3}+\cfrac{n-22}{d_{4}+Mn}}}\label{eq:2.9}
\end{equation}

\par\end{flushleft}

with $d_{i},M\in\mathbb{R}$ and $i,n\in\mathbb{Z}^{+}$, $1\leq i\leq4$.
The coefficients $d_{i}$ and $M$ are given in Appendix 3. For $n=25$,
(\ref{eq:2.9}) yields $\left|\psi_{26}\right|\approx7.07767\times10^{-9}$.
Assuming $\psi_{26}\approx\pm7.1\times10^{-9}$ yields a conservative
assumption for the value of $\pi(10^{26})$
\begin{equation}
\pi(10^{26})\approx1699246738822618041025224\pm12064651845640588\label{eq:2.10}
\end{equation}

or
\begin{equation}
1699246726757966195384636\lesssim\pi(10^{26})\lesssim1699246750887269886665812\label{eq:2.11}
\end{equation}

Assuming further that $\psi_{26}$ would be positive like $\delta_{25}^{\lyxmathsym{\textasciiacute}}$
and $\delta_{24}^{\lyxmathsym{\textasciiacute}}$ and that it would
be included in the range $7\times10^{-9}\lesssim\psi_{26}\lesssim7.1\times10^{-9}$,
the value range for $\pi(10^{26})$ reduces to 
\begin{equation}
\pi(10^{26})\approx1699246738822618041025224_{+11894727171758326}^{+12064651845640588}\label{eq:2.12}
\end{equation}

or 
\begin{equation}
1699246750717345212783550\lesssim\pi(10^{26})\lesssim1699246750887269886665812\label{eq:2.13}
\end{equation}

\section{Discussion}

It is interesting to compare this value range to those values obtained
with the approximating functions of Section 1. Table 2 shows for $24\leq n\leq26$
the values of $Round\left(\frac{10^{n}}{log\left(10^{n}\right)}\right)$,
$Round\left(Li(10^{n})\right)$ (\ref{eq:1.2}), $Round\left(R\left(10^{n}\right)\right)$
(\ref{eq:1.3}) and the associated relative differences 
\begin{equation}
\delta_{n}^{\lyxmathsym{\textacutedbl}}=\frac{\pi(10^{n})-Round\left(f(n)\right)}{\pi(10^{n})}\label{eq:3.8}
\end{equation}

where $f(n)$ is either $\left(\frac{10^{n}}{log\left(10^{n}\right)}\right)$,
or $Li(10^{n})$, or $R\left(10^{n}\right)$, and where $\pi(10^{24})=18435599767349200867866$
and $\pi(10^{25})=176846309399143769411680$.

\begin{table}
\caption{Values of $Round\left(\frac{10^{n}}{log\left(10^{n}\right)}\right)$,
$Round\left(Li(10^{n})\right)$ (\ref{eq:1.2}), $Round\left(R\left(10^{n}\right)\right)$
(\ref{eq:1.3}) for $24\leq n\leq26$ and associated $\delta_{n}^{\lyxmathsym{\textacutedbl}}$
(\ref{eq:2.12})}

\begin{tabular}{|c|c|c|c|}
\hline 
$n$ & $f(n)$ & $Round\left(f(n)\right)$ & $\delta_{n}^{\lyxmathsym{\textacutedbl}}$\tabularnewline
\hline 
\hline 
24 & $\frac{10^{n}}{log\left(10^{n}\right)}$ & 18095603412635492818797 & $1.84424\times10^{-2}$\tabularnewline
\hline 
24 & $Li(10^{n})$ & 18435599767366347775143 & $-9.30098\times10^{-13}$\tabularnewline
\hline 
24 & $R\left(10^{n}\right)$ & 18435599767347541878147 & $8.99884\times10^{-14}$\tabularnewline
\hline 
\hline 
25 & $\frac{10^{n}}{log\left(10^{n}\right)}$ & 173717792761300731060452 & $1.76906\times10^{-2}$\tabularnewline
\hline 
25 & $Li(10^{n})$ & 176846309399198930392618 & $-3.11915\times10^{-13}$\tabularnewline
\hline 
25 & $R\left(10^{n}\right)$ & 176846309399141934626966 & $1.03750\times10^{-14}$\tabularnewline
\hline 
\hline 
26 & $\frac{10^{n}}{log\left(10^{n}\right)}$ & 1670363391935583952504342 & -\tabularnewline
\hline 
26 & $Li(10^{n})$ & 1699246750872593033005723 & -\tabularnewline
\hline 
26 & $R\left(10^{n}\right)$ & 1699246750872419991992147 & -\tabularnewline
\hline 
\end{tabular}
\end{table}

Note that in calculating $R\left(10^{n}\right)$, the summation in
(\ref{eq:1.3}) is made for $1\leq j\leq1000$, as the value of $R\left(10^{n}\right)$
does not change for higher values of $j$, due mainly to the operation
of rounding to the nearest integer. Similarly, the values obtained
by rounding to the nearest integer of the function (\ref{eq:1.4})
does not differ from those obtained by rounding $R\left(10^{n}\right)$
for the same reason.

It is seen that, for $n=24$ and $25$, the values obtained with $\frac{10^{n}}{log\left(10^{n}\right)}$
are relatively far from the real values of $\pi(10^{n})$ with $\delta_{n}^{\lyxmathsym{\textacutedbl}}$
of the order $10^{-2}$ , while the values of $Li(10^{n})$ and $R\left(10^{n}\right)$
gives much better results with $\delta_{n}^{\lyxmathsym{\textacutedbl}}$
of the order $10^{-13}...\,10^{-14}$. For $n=26$, the values of
$Li(10^{26})$ and $R\left(10^{26}\right)$ are within the range given
by (\ref{eq:2.13}).

\section{Conclusion}

We proposed a simple method to interpolate a conjectured value of
$\pi\left(10^{26}\right)$, based on polynomial interpolations of
the known 25 first values of $\pi\left(10^{n}\right)$ with two corrective
functions, the first based on a Thiele interpolation of the relative
differences of the exact and interpolated values of $\pi\left(10^{n}\right)$
and the second obtained by a new Thiele interpolation of the last
five values of the relative differences of the exact and interpolated
values of $\pi\left(10^{n}\right)$ with a previous value of the first
corrective function. The range obtained for the conjectured value
of $\pi\left(10^{26}\right)$ contains those values calculated by
the Eulerian logarithmic integral and with the Riemann functions.

\section*{Appendix 1: Polynomials $P_{n}\left(x\right)$ for $2\leq n\leq25$}

$P_{2}\left(x\right)=21x-17$

\medskip{}
$P_{3}\left(x\right)=61x^{2}-162x+105$\medskip{}

$P_{4}\left(x\right)=\left(398/3\right)x^{3}-735x^{2}+\left(3892/3\right)x-691$\medskip{}

$P_{5}\left(x\right)=\left(1397/6\right)x^{4}-\left(6587/3\right)x^{3}+\left(44485/6\right)x^{2}-\left(31033/3\right)x+4897$\medskip{}

$P_{6}\left(x\right)=\left(41269/120\right)x^{5}-\left(118219/24\right)x^{4}+\left(648877/24\right)x^{3}-\left(1679165/24\right)x^{2}$

$+\left(1677731/20\right)x-36372$\medskip{}

$P_{7}\left(x\right)=\left(63053/144\right)x^{6}-\left(2124317/240\right)x^{5}+\left(10324961/144\right)x^{4}-\left(14150231/48\right)x^{3}$

$+\left(46161541/72\right)x^{2}-\left(6885127/10\right)x+278893$\medskip{}

$P_{8}\left(x\right)=\left(1230899/2520\right)x^{7}-\left(1059103/80\right)x^{6}+\left(106869757/720\right)x^{5}-\left(42511909/48\right)x^{4}$

$+\left(2168305201/720\right)x^{3}-\left(692786273/120\right)x^{2}+\left(199310224/35\right)x-2182905$\medskip{}

$P_{9}\left(x\right)=\left(1087519/2240\right)x^{8}-\left(12232463/720\right)x^{7}+\left(8059007/32\right)x^{6}$

$-\left(295746589/144\right)x^{5}+\left(9612782119/960\right)x^{4}-\left(5337867053/180\right)x^{3}$

$+\left(17329512589/336\right)x^{2}-\left(190033621/4\right)x+17392437$\medskip{}

$P_{10}\left(x\right)=\left(8766223/20160\right)x^{9}-\left(1526557/80\right)x^{8}+\left(3642052523/10080\right)x^{7}$

$-\left(61717169/16\right)x^{6}+\left(73322957917/2880\right)x^{5}-\left(6425876741/60\right)x^{4}$

$+\left(119710398973/420\right)x^{3}-\left(11000451139/24\right)x^{2}+\left(13960806758/35\right)x-140399577$\medskip{}

$P_{11}\left(x\right)=\left(1285948381/3628800\right)x^{10}-\left(13829848163/725760\right)x^{9}+\left(2713678729/6048\right)x^{8}$

$-\left(734294140229/120960\right)x^{7}+\left(8994784761253/172800\right)x^{6}-\left(10167707046167/34560\right)x^{5}$

$+\left(80107171523893/72576\right)x^{4}-\left(488994253144639/181440\right)x^{3}+\left(204682947179773/50400\right)x^{2}$

$-\left(1697281382717/504\right)x+1145548804$\medskip{}

$P_{12}\left(x\right)=\left(1762482209/6652800\right)x^{11}-\left(62163411143/3628800\right)x^{10}+\left(356291415727/725760\right)x^{9}$

$-\left(496320428783/60480\right)x^{8}+\left(1914173704139/21600\right)x^{7}-\left(111749346412919/172800\right)x^{6}$

$+\left(261253448575099/80640\right)x^{5}-\left(4021281732287857/362880\right)x^{4}+\left(22852595048399971/907200\right)x^{3}$

$-\left(1810394687049689/50400\right)x^{2}+\left(113126344733737/3960\right)x-9429344450$\medskip{}

$P_{13}\left(x\right)=\left(87613503863/479001600\right)x^{12}-\left(1117825763711/79833600\right)x^{11}$

$+\left(4178914904089/8709120\right)x^{10}-\left(671147586733/69120\right)x^{9}$

$+\left(1870673383324673/14515200\right)x^{8}-\left(2850595750775761/2419200\right)x^{7}$

$+\left(66036047077439107/8709120\right)x^{6}-\left(50008128480417883/1451520\right)x^{5}$

$+\left(1188002313137637187/10886400\right)x^{4}-\left(2241074007266561/9600\right)x^{3}$

$+\left(105572714697948823/332640\right)x^{2}-\left(562059733555247/2310\right)x+78184159413$\medskip{}

$P_{14}\left(x\right)=\left(730505748461/6227020800\right)x^{13}-\left(418827227947/39916800\right)x^{12}$

$+\left(202948195226041/479001600\right)x^{11}-\left(37036465604101/3628800\right)x^{10}$

$+\left(338543894892361/2073600\right)x^{9}-\left(2209549040385367/1209600\right)x^{8}$

$+\left(638513890320991673/43545600\right)x^{7}-\left(309191577104894683/3628800\right)x^{6}$

$+\left(3885522866329964983/10886400\right)x^{5}-\left(1921510055977553161/1814400\right)x^{4}$

$+\left(1532239610437270573/712800\right)x^{3}-\left(517946743840607861/184800\right)x^{2}$

$+\left(749473155761448241/360360\right)x-652321589048$\medskip{}

$P_{15}\left(x\right)=\left(68044412311/968647680\right)x^{14}-\left(434610313091/59875200\right)x^{13}$

$+\left(54461331244787/159667200\right)x^{12}-\left(1154924881829021/119750400\right)x^{11}$

$+\left(2655624146858351/14515200\right)x^{10}-\left(1274125884612191/518400\right)x^{9}$

$+\left(814648425163663999/33868800\right)x^{8}-\left(1890972821539241113/10886400\right)x^{7}$

$+\left(3363912482966950097/3628800\right)x^{6}-\left(39446931450255536567/10886400\right)x^{5}$

$+\left(202795315040803498879/19958400\right)x^{4}-\left(28079184621772453729/1425600\right)x^{3}$

$+\left(416132103978867941789/16816800\right)x^{2}-\left(70617762061263719/3960\right)x$

$+5471675518942$\medskip{}

$P_{16}\left(x\right)=\left(51581436427121/1307674368000\right)x^{15}-\left(29037678164927/6227020800\right)x^{14}$

$+\left(2356528303232491/9340531200\right)x^{13}-\left(3963130920435319/479001600\right)x^{12}$

$+\left(188475479592111883/1026432000\right)x^{11}-\left(42323388417591161/14515200\right)x^{10}$

$+\left(15614861323300656127/457228800\right)x^{9}-\left(13050538335127845887/43545600\right)x^{8}$

$+\left(2588220260396190960353/1306368000\right)x^{7}-\left(53527147712466477239/5443200\right)x^{6}$

$+\left(232250648536987989971/6415200\right)x^{5}-\left(214241913350759156389/2217600\right)x^{4}$

$+\left(544191268377690363241837/3027024000\right)x^{3}-\left(524844387326399456539/2402400\right)x^{2}$

$+\left(6906580916276240896/45045\right)x-46109760908179$\medskip{}

$P_{17}\left(x\right)=\left(17450310541093/836911595520\right)x^{16}-\left(7313219107110283/2615348736000\right)x^{15}$

$+\left(180531219567405389/1046139494400\right)x^{14}-\left(121365437052931943/18681062400\right)x^{13}$

$+\left(3833898578147855953/22992076800\right)x^{12}-\left(44388621420665486663/14370048000\right)x^{11}$

$+\left(311968634799119772881/7315660800\right)x^{10}-\left(408405984743715765221/914457600\right)x^{9}$

$+\left(104689503078054721538951/29262643200\right)x^{8}-\left(57304152503501292804619/2612736000\right)x^{7}$

$+\left(587617527762868915556291/5748019200\right)x^{6}-\left(256560216571537768707173/718502400\right)x^{5}$

$+\left(981414623224315552474067/1076275200\right)x^{4}-\left(29710845425841425154523003/18162144000\right)x^{3}$

$+\left(389276944856756983375399/201801600\right)x^{2}-\left(105829165042439513371/80080\right)x$

$+390148002619146$\medskip{}

$P_{18}\left(x\right)=\left(1851587512783157/177843714048000\right)x^{17}-\left(32892317466569501/20922789888000\right)x^{16}$

$+\left(256329638772689/2335132800\right)x^{15}-\left(4918740790637408989/1046139494400\right)x^{14}$

$+\left(103257725801778476353/747242496000\right)x^{13}-\left(339479305160332664821/114960384000\right)x^{12}$

$+\left(636260413866908181391/13412044800\right)x^{11}-\left(170936402957313265801/292626432\right)x^{10}$

$+\left(58215415330924617167453/10450944000\right)x^{9}-\left(6039147285766985737916771/146313216000\right)x^{8}$

$+\left(1362372951548947987220519/5748019200\right)x^{7}-\left(6008945224011701433565219/5748019200\right)x^{6}$

$+\left(2278512605057461714112724299/653837184000\right)x^{5}$

$-\left(1243038724423425067320005221/145297152000\right)x^{4}$

$+\left(2568389927503727128176139/172972800\right)x^{3}-\left(3437947478393805445387967/201801600\right)x^{2}$

$+\left(139868165483953060330259/12252240\right)x-3313027022947168$\medskip{}

$P_{19}\left(x\right)=\left(2628028871421403/533531142144000\right)x^{18}-\left(10567575582731201/12703122432000\right)x^{17}$

$+\left(21003369878917183/321889075200\right)x^{16}-\left(8251376607822726667/2615348736000\right)x^{15}$

$+\left(7058709267627854321/67060224000\right)x^{14}-\left(137168476397716450849/53374464000\right)x^{13}$

$+\left(22944546189331191060451/482833612800\right)x^{12}-\left(271196454331659216727/399168000\right)x^{11}$

$+\left(556013030981016791428933/73156608000\right)x^{10}-\left(2446752331753149321520481/36578304000\right)x^{9}$

$+\left(1840753991498633195715463/3973939200\right)x^{8}-\left(72170805307845647867807051/28740096000\right)x^{7}$

$+\left(82740150719665463161579783529/7846046208000\right)x^{6}$

$-\left(98402301347766134600009057/2918916000\right)x^{5}$

$+\left(107115515185413881741729737/1341204480\right)x^{4}$

$-\left(814304007652526714435690759/6054048000\right)x^{3}$

$+\left(119147362970798190946403539/791683200\right)x^{2}-\left(1210609739820985493184649/12252240\right)x$

$+28223319434109668$\medskip{}

$P_{20}\left(x\right)=\left(26934382640709253/12164510040883200\right)x^{19}$

$-\left(166368869850866779/400148356608000\right)x^{18}$

$+\left(2987733695853594407/82081714176000\right)x^{17}$

$-\left(3538361726199431459/1793381990400\right)x^{16}$

$+\left(211578849172220030591/2853107712000\right)x^{15}$

$-\left(32187947021596169826761/15692092416000\right)x^{14}$

$+\left(313223082576559198263053/7242504192000\right)x^{13}$

$-\left(1028851538290603359043297/1448500838400\right)x^{12}$

$+\left(44456420445740629307608553/4828336128000\right)x^{11}$

$-\left(212467505628141678941887/2239488000\right)x^{10}$

$+\left(3760078141213234941907434283/4828336128000\right)x^{9}$

$-\left(4906343664090048465680852441/965667225600\right)x^{8}$

$+\left(4812537157115295300593562547/183891708000\right)x^{7}$

$-\left(2475493108822130612914163579363/23538138624000\right)x^{6}$

$+\left(986767217403394741887962341/3048192000\right)x^{5}$

$-\left(5547876667800341701953704701/7472424960\right)x^{4}$

$+\left(125190395324618857278613349737/102918816000\right)x^{3}$

$-\left(41062191060690947895024739489/30875644800\right)x^{2}$

$+\left(15342006664059599043795083/17907120\right)x-241120506972982862$\medskip{}

$P_{21}\left(x\right)=\left(768801874165960931/810967336058880000\right)x^{20}$

$-\left(47895830419641353593/243290200817664000\right)x^{19}$

$+\left(19593696496463403689/1024379792916480\right)x^{18}$

$-\left(448205658674858604491/388022648832000\right)x^{17}$

$+\left(8712644329419234113383/179338199040000\right)x^{16}$

$-\left(189707965396142126979173/125536739328000\right)x^{15}$

$+\left(1809328805921664964890847/50214695731200\right)x^{14}$

$-\left(253665818771705300538215719/376610217984000\right)x^{13}$

$+\left(527355591384145604585438497/52672757760000\right)x^{12}$

$-\left(460925183825799784732199789/3862668902400\right)x^{11}$

$+\left(252657778262904116233611461/220723937280\right)x^{10}$

$-\left(170656737803291291159227967873/19313344512000\right)x^{9}$

$+\left(1906512687609116773290392565521/34871316480000\right)x^{8}$

$-\left(2302249554842995226554178587199/8559323136000\right)x^{7}$

$+\left(19560002409333744035028716869169/18830510899200\right)x^{6}$

$-\left(199482007364920203714291611479/64576512000\right)x^{5}$

$+\left(54603156665489883803994128245399/7939451520000\right)x^{4}$

$-\left(753730761427013331958837118357/68612544000\right)x^{3}$

$+\left(156740472492042038497972499701/13332664800\right)x^{2}$

$-\left(866115476410575028425322373/116396280\right)x+2065285115524899931$\medskip{}

$P_{22}\left(x\right)=\left(54386159471042293/140359731240960000\right)x^{21}$

$-\left(11339777731555550441/128047474114560000\right)x^{20}$

$+\left(6933583440707809527311/729870602452992000\right)x^{19}$

$-\left(16277845641786522178019/25609494822912000\right)x^{18}$

$+\left(173069291634903343162099/5820339732480000\right)x^{17}$

$-\left(185313460297730292870869/179338199040000\right)x^{16}$

$+\left(31253190030180758167460353/1129830653952000\right)x^{15}$

$-\left(439746760244290352119569907/753220435968000\right)x^{14}$

$+\left(27854308762263416011926122669/2824576634880000\right)x^{13}$

$-\left(7087441471779667149333596851/52672757760000\right)x^{12}$

$+\left(17272738385198816463545470907/11588006707200\right)x^{11}$

$-\left(74078609744831895304816687447/5518098432000\right)x^{10}$

$+\left(85271649858214871482393447550287/869100503040000\right)x^{9}$

$-\left(544996898468027472985563167444741/941525544960000\right)x^{8}$

$+\left(70149627027741059244452076836923/25677969408000\right)x^{7}$

$-\left(957972246454320768062382326886223/94152554496000\right)x^{6}$

$+\left(39105801276089863543963319222962243/1333827855360000\right)x^{5}$

$-\left(504262702012069575930635459206387/7939451520000\right)x^{4}$

$+\left(71182532364843711432355416393311/718331328000\right)x^{3}$

$-\left(36474259397560749064143053993/350859600\right)x^{2}$

$+\left(7533688390536925014719660137/116396280\right)x-17731276931934494721$\medskip{}

$P_{23}\left(x\right)=\left(626047699915867799/4132355616829440000\right)x^{22}$

$-\left(32308010721489584677/851515702861824000\right)x^{21}$

$+\left(2511103473848270746499/561438924963840000\right)x^{20}$

$-\left(990899107137455767879/3003582726144000\right)x^{19}$

$+\left(2186295775571121711474881/128047474114560000\right)x^{18}$

$-\left(255701609845973121847979/388022648832000\right)x^{17}$

$+\left(119722559289831616224557863/6083703521280000\right)x^{16}$

$-\left(1227534164901745056978902207/2636271525888000\right)x^{15}$

$+\left(100182659337689089936042566167/11298306539520000\right)x^{14}$

$-\left(2154242664974814138465530011/15692092416000\right)x^{13}$

$+\left(3693244759563863089340159781247/2124467896320000\right)x^{12}$

$-\left(348679419716275796227603114639/19313344512000\right)x^{11}$

$+\left(24313863183439316679673602226004513/158176291553280000\right)x^{10}$

$-\left(225569651273130061058657092854493/210901722071040\right)x^{9}$

$+\left(187595874338792344136694981328643/31039303680000\right)x^{8}$

$-\left(3134444555753150961218198128157/114124308480\right)x^{7}$

$+\left(46582581480572090585766685040140853/470762772480000\right)x^{6}$

$-\left(24628252584588755381459107318208291/88921857024000\right)x^{5}$

$+\left(332629942234296028174359960597644779/568586874240000\right)x^{4}$

$-\left(73343429897691545576959449692313389/82129215168000\right)x^{3}$

$+\left(741767456331546132245359568800763/806626220400\right)x^{2}$

$-\left(65620214252713224536332366223/116396280\right)x+152553697445181546607$\medskip{}

$P_{24}\left(x\right)=\left(97841530496558649293/1723467782592331776000\right)x^{23}$

$-\left(335407509903912323681/21615398611107840000\right)x^{22}$

$+\left(454501108157980384073/227070854096486400\right)x^{21}$

$-\left(107504314873384627683143/663518729502720000\right)x^{20}$

$+\left(8983776081040688437620797/973160803270656000\right)x^{19}$

$-\left(50404505281634117870005549/128047474114560000\right)x^{18}$

$+\left(1168666259954948513855852177/89633231880192000\right)x^{17}$

$-\left(2471024795530257064602656791/7189831434240000\right)x^{16}$

$+\left(51511683265529304575508645967/7030057402368000\right)x^{15}$

$-\left(721535570132495990554385229809/5649153269760000\right)x^{14}$

$+\left(30370138820999258718781973500777/16570849591296000\right)x^{13}$

$-\left(46160207442983253249607104925813/2124467896320000\right)x^{12}$

$+\left(408425601328540276858577344312109/1917288382464000\right)x^{11}$

$-\left(1907340749270059838043121409492039/1106127912960000\right)x^{10}$

$+\left(10084417668013132593085446666361763/878757175296000\right)x^{9}$

$-\left(898532173217897912413387346270287/14411105280000\right)x^{8}$

$+\left(876114367621403488143688018877125577/3201186852864000\right)x^{7}$

$-\left(695775343904264963827651396107198739/727542466560000\right)x^{6}$

$+\left(26751505111131816927397223600866889/10266151896000\right)x^{5}$

$-\left(79477104547597150222616625089859530191/14783258730240000\right)x^{4}$

$+\left(403589896543329394729142860419921947/50190075936000\right)x^{3}$

$-\left(52511426579628409245398898318755171/6453009763200\right)x^{2}$

$+\left(26325492596728861363119710723467/5354228880\right)x-1315069260003198192788$\medskip{}

$P_{25}\left(x\right)=\left(1810237126686092753419/88635485961891348480000\right)x^{24}$

$-\left(62771250251033894473907/10340806695553990656000\right)x^{23}$

$+\left(11511991109112409288411081/13488008733331292160000\right)x^{22}$

$-\left(308610477223222584935011/4087275373736755200\right)x^{21}$

$+\left(1651947699182756157926795881/350337889177436160000\right)x^{20}$

$-\left(430292416333562341190138281/1946321606541312000\right)x^{19}$

$+\left(1693062186931273696473576313/209532230369280000\right)x^{18}$

$-\left(2440242314982415470323718869/10342295986176000\right)x^{17}$

$+\left(42496596088997036408436172662149/7592461994557440000\right)x^{16}$

$-\left(13805367401638213159287611947619/126541033242624000\right)x^{15}$

$+\left(10501406884609102136813341938733721/5965505852866560000\right)x^{14}$

$-\left(260923334002982532857433489649357/11047233060864000\right)x^{13}$

$+\left(3149961803456717214034882917658324207/11931011705733120000\right)x^{12}$

$-\left(311007895030981777959409995030791881/126541033242624000\right)x^{11}$

$+\left(36072334629973139334722245407290958851/1898115498639360000\right)x^{10}$

$-\left(3842017793050978223673061999182697157/31635258310656000\right)x^{9}$

$+\left(9526409805672264569123818314621843047/14966587883520000\right)x^{8}$

$-\left(17325638483889432743363938986730301171/6402373705728000\right)x^{7}$

$+\left(5590963217656815416105678129476360121221/608225502044160000\right)x^{6}$

$-\left(4938656200253742430432451463654529907/202164221952000\right)x^{5}$

$+\left(16036374727879322218852223186722596976589/325231692065280000\right)x^{4}$

$-\left(43579059908227428615914106531732753781/602280911232000\right)x^{3}$

$+\left(3563764005018304979715944223158462749/49473074851200\right)x^{2}$

$-\left(114930909252025411165097498628379/2677114440\right)x+11356590626799451081145$

\section*{Appendix 2: Coefficients $c_{i}$ and $K$ of the corrective function
$\Phi_{25}$ }

$c_{1}=0.72619047619047619047619047619047619047619047619048$

$c_{2}=-12.737322640345465761875385564466378778531770512029$

$c_{3}=-0.84164269999731510129175214910348073427241614034768$

$c_{4}=-17.908330903689032603895485760334710528008413393816$

$c_{5}=-21.793051968869308781494503792492773700351128057503$

$c_{6}=0.052225391819545299875787414478370692540986477002$

$c_{7}=22.867256875351930467366237953412010183478716371242$

$c_{8}=-1.707611319924344890707228312842658677594869606318$

$c_{9}=-1.2751632549654008209020770141453274057116453048774$

$c_{10}=92.538545141732012060260360179516606702270761634364$

$c_{11}=0.15088559448367433870497668805212911805725048971437$

$c_{12}=12.296035579537995708862927240169729240327162730035$

$c_{13}=0.47119594517891782795780492658781237198369447024440$

$c_{14}=103.88603190326434306252340497028079649950656584049$

$c_{15}=0.00206293819355223802469982625483902305678068096840$

$c_{16}=-187.13045591233740512137368315625509082987760574938$

$c_{17}=-0.27349669744110325871269828279318826062620906588370$

$c_{18}=60.84802398138390071856628905328286526305777316700$

$c_{19}=0.038078160763329325175635232207793882248458501401699$

$c_{20}=-375.10551758837650646142692040134675795455506670138$

$c_{21}=0.12511512946150416165180535229285204189299499111780$

$c_{22}=71.229908224807213192533908002847675276237637077608$

$K=-2.4774348845851811138845060795239470850970682780195$

\section*{Appendix 3: Coefficients $d_{i}$ and $M$ of the corrective function
$\left|\psi_{n+1}\right|$}

$d_{1}=-3.6415421544257471159022604244391523071564117009437\times10^{-7}$

$d_{2}=-3485608.6856471123162187714001336338709185526509134$

$d_{3}=-4.6018631805033582153977197516014111986555627816497\times10^{-7}$

$d_{4}=-294118428.93353607774712567420242866941568199754406$

$M=13267560.424903927928414820656294875945590431439845$


\begin{thebibliography}{References}
\bibitem{1} D. Zagier (1997). \textquotedbl{}Newman's short proof
of the prime number theorem\textquotedbl{}. American Mathematical
Monthly 104 (8): 705\textendash{}708. doi:10.2307/2975232.

\bibitem[2]{2} E.W. Weisstein. \textquotedbl{}Prime Counting Function.\textquotedbl{}
From MathWorld--A Wolfram Web Resource. http://mathworld.wolfram.com/PrimeCountingFunction.html,
last accessed 25/02/2013.

\bibitem[3]{3} E.W. Weisstein. \textquotedbl{}Riemann Prime Counting
Function.\textquotedbl{} From MathWorld--A Wolfram Web Resource. http://mathworld.wolfram.com/RiemannPrimeCountingFunction.html,
last accessed 25/02/2013.

\bibitem[4]{4} Séroul, R. \textquotedbl{}Legendre's Formula\textquotedbl{}
and \textquotedbl{}Implementation of Legendre's Formula.\textquotedbl{}
§8.7.1 and 8.7.2 in Programming for Mathematicians. Berlin: Springer-Verlag,
pp. 175-179, 2000.

\bibitem[5]{5} Riesel, H. \textquotedbl{}Lehmer's Formula.\textquotedbl{}
Prime Numbers and Computer Methods for Factorization, 2nd ed. Boston,
MA: Birkhäuser, pp. 13-14, 1994.

\bibitem[6]{6} Séroul, R. \textquotedbl{}Meissel's Formula.\textquotedbl{}
§8.7.3 in Programming for Mathematicians. Berlin: Springer-Verlag,
pp. 179-181, 2000.

\bibitem[7]{7} N. J. A. Sloane, S. Plouffe, \textquotedbl{}Number
of primes < 10\textasciicircum{}n\textquotedbl{}, Sequence A006880
in \emph{The On-line Encyclopedia of Integer Sequences}, published
electronically at http://oeis.org, last accessed 7 June 2013.

\bibitem[8]{8} Jens Franke et al., 29 July 2010, in \cite{7}

\bibitem[9]{9} J. Buethe, J. Franke, A. Jost, T. Kleinjung, 01 June
2013, in \cite{7}\end{thebibliography}
\end{document}